\documentclass[10pt,reqno]{amsart}
\usepackage[T2A]{fontenc}
\usepackage{amsmath}
\usepackage{amssymb}
\usepackage{amsfonts}
\usepackage{bm} % for \bm
\usepackage{stmaryrd} % for \trianglelefteqslant

\newtheorem{thm}{Theorem}
\newtheorem{prop}[thm]{Proposition}
\newtheorem{lem}[thm]{Lemma}
\newtheorem{cor}[thm]{Corollary}

\begin{document}

\title[Automorphisms of nonsplit coverings]
{Automorphisms of nonsplit coverings of $PSL_2(q)$ \\ in odd characteristic dividing $q-1$}

\author{Andrei V. Zavarnitsine}%
\address{Andrei V. Zavarnitsine
\newline\indent Sobolev Institute of Mathematics,
\newline\indent 4, Koptyug av.
\newline\indent 630090, Novosibirsk, Russia
} \email{zav@math.nsc.ru}

\maketitle {\small
\begin{quote}
\noindent{\sc Abstract. } We classify the nonsplit extensions of elementary abe\-lian $p$-groups
by $\operatorname{PSL}_2(q)$, with odd $p$ dividing $q-1$, for an irreducible induced action,
calculate the relevant low-dimensional cohomology groups, and describe the automorphism groups of such extensions.
\medskip

 \end{quote}
}

\section{Introduction}

Given a short exact sequence of groups
\begin{equation}\label{g}
0\to V \to G \to L\to 1,
\end{equation}
where $V$ is abelian (written additively), we say that $G$ is {\em an extension} of $V$ by~$L$, or
{\em a covering} of $L$ with kernel $V$.
Such extensions arise naturally in inductive arguments
or when constructing minimal examples and counterexamples. We will be interested in the case
where $G$ is finite and nonsplit and $V$ acquires the structure of an irreducible $FL$-module (for a
suitable finite field $F$ of characteristic $p$) from the conjugation in $G$.
Such extensions can only exist if $p$ divides $|L|$.
We also restrict ourselves to the case $L\cong\operatorname{PSL}_2(q)$.
Extensions of this form for $p=2$ and $q$ odd were explicitly constructed in \cite{00Buri}, and their automorphism groups were described \cite{21RevZav}. Some results in the case of $q$ being a power of $p$ were obtained in~\cite{15Bur}.
The aim of this paper is to classify such extensions in the case $2\ne p\mid (q-1)$ and describe their automorphism groups.
In this case, we can use the fact that the natural permutation $FL$-module arising from the action of $L$ on the projective
line over $\mathbb{F}_q$ is completely reducible. This is not so if $2\ne p\mid (q+1)$
which case will be a subject of future research. We now state the main results.

\begin{thm} \label{nonsplit} Up to isomorphism there is a unique nonsplit extension
of an ele\-men\-ta\-ry abelian $p$-group~$V$ by
$L=\operatorname{PSL}_2(q)$ with irreducible induced action of $L$ on~$V$, where $2\ne p\mid (q-1)$.
In this extension, $|V|=p^q$.
\end{thm}

The group $V$ from Theorem \ref{nonsplit} as an $\mathbb{F}_pL$-module can be identified with the unique
nonprincipal irreducible module in the principal $p$-block of $L$.
The low-dimensional co\-ho\-mo\-lo\-gy of $V$ is as follows.

\begin{thm}\label{2coh}
In the above notation, we have
$H^1(L,V)\cong H^2(L,V)\cong \mathbb{F}_p$.
\end{thm}

Recall that $P\Gamma L_2(q)$ denotes the extension of $PGL_2(q)$ by its field auto\-mor\-phisms.
The automorphism group of the nonsplit extension from Theorem \ref{nonsplit} is described~by

\begin{thm} \label{aut}
Let $G$ fit in the nonsplit exact sequence $(\ref{g})$, where $V$ is an irreducible $\mathbb{F}_pL$-module
for $L=\operatorname{PSL}_2(q)$ and $2\ne p\mid (q-1)$.  Then there is a short exact sequence
\begin{equation}\label{autg}
\qquad 0\to W \to \operatorname{Aut}(G) \to P\Gamma L_2(q)\to 1,
\end{equation}
where $W$ is elementary abelian of order $p^{q+1}$.
\end{thm}

\section{Auxiliary facts}

Basic notation and facts of homological algebra can be found in \cite{70Gru,94Wei}. For abelian groups $A$ and $B$, we denote $\operatorname{Hom}(A,B) = \operatorname{Hom}_\mathbb{Z}(A,B)$ and $\operatorname{Ext}(A,B) = \operatorname{Ext}_\mathbb{Z}^1(A,B)$.

\begin{lem}[The Universal Coefficient Theorem for Cohomology]\cite[Ch.\,3, The\-o\-rem 3]{70Gru} \label{ucoe} For all $i\geqslant 1$, every group $G$, and every trivial $G$-module $A$,
$$
H^i(G,A)\cong\operatorname{Hom}(H_i(G,\mathbb{Z}),A)\oplus\operatorname{Ext}(H_{i-1}(G,\mathbb{Z}),A).
$$
\end{lem}
\begin{lem}\cite[\S 3.5]{70Gru}\label{coh1} For a trivial $G$-module $A$, we have
\begin{enumerate}
  \item[$(i)$] $H^1(G,A)\cong \operatorname{Hom}(G/G',A)$.
  \item[$(ii)$] $H_1(G,A)\cong G/G'\otimes_\mathbb{Z} A$.
\end{enumerate}
\end{lem}

\begin{lem}[Shapiro's lemma]\cite[\S 6.3]{94Wei} \label{shap} Let $H\leqslant G$ with $|G:H|$ finite. If $V$ is an $H$-module and $i\geqslant 0$ then $H^i(G,V^G)\cong H^i(H,V)$, where $V^G$ is the induced $G$-module.
\end{lem}

\begin{lem}\cite[p. 322]{18Deo}\label{ezmn} $\operatorname{Ext}(\mathbb{Z}_m,\mathbb{Z}_n)\cong \mathbb{Z}_d$, where $d=(m,n)$.
\end{lem}

\begin{lem}\cite[Proposition 3.3.4]{94Wei}\label{esum}.
$\operatorname{Ext}_R^i(A, B_1\oplus B_2)\cong \operatorname{Ext}_R^i(A,B_1)\oplus \operatorname{Ext}_R^i(A,B_2)$ for all rings $R$, $R$-modules $A,B_1,B_2$, and all $i\geqslant 0$.
\end{lem}

The Schur multiplier of a group $G$ is denoted by $\operatorname{Sch}(G)$.
If $A$ is a finite abelian group and $p$ a prime then $A_{(p)}$ denotes the $p$-primary component of $A$.
Henceforth, we assume that $G$ is finite.

\begin{lem}\cite[Theorem 25.1]{67Hup} \label{tim} Let $p$ be a prime and let $S\in\operatorname{Syl}_p(G)$.
Then $\operatorname{Sch}(G)_{(p)}$ is isomorphic to a subgroup of $\operatorname{Sch}(S)$.
\end{lem}

\begin{lem}\cite{76Gag}\label{lgag}
Let $F$ be a field of characteristic $p>0$ and let $V$ be an irreducible
$FG$-module that does not belong to the principal $p$-block of $G$. Then $H^n(G,V) = 0$ for all $n\geqslant 0$.
\end{lem}

Let $\theta$ be an irreducible character of $G$.
If $\operatorname{Z}(G)=1$ then $G\trianglelefteqslant \operatorname{Aut}(G)$
and we may speak of the inertia group
$I_{\operatorname{Aut}(G)}(\theta)=\{g\in \operatorname{Aut}(G)\mid \theta^g=\theta\}$.

\begin{prop}\cite[Proposition 4]{21RevZav}\label{normim}
Let $F$ be a field and $\mathcal{X}$ a faithful irreducible $F$-representation of a group $G$
with Brauer character $\theta\in\operatorname{iBr}_F(G)$ of degree $n$. Suppose that
$\operatorname{Z}(G)=1$ and denote
$$N=N_{\operatorname{GL}_n(F)}(\mathcal{X}(G))\quad \text{and} \quad Z=C_{\operatorname{GL}_n(F)}(\mathcal{X}(G)).$$ Then
$N/Z\cong I_{\operatorname{Aut}(G)}(\theta).$
\end{prop}

\section{Isomorphic extensions}

Let $Q$ be a group, $K$ a commutative ring with $1$, and $M$ a right $KQ$-module. The pair
$(\nu,\mu)\in \operatorname{Aut}(Q)\times \operatorname{Aut}_K(M)$ is {\it compatible} if
$$
(mg)\mu = (m\mu)(g\nu)
$$
for all $m\in M$, $g\in Q$. The set of all compatible pairs forms a group $\operatorname{Comp}(Q,M)$ under composition. Given $\tau\in Z^2(Q,M)$, one can define
\begin{equation}\label{cact}
\tau^{(\nu,\mu)}(g,h)=\tau(g\nu^{-1},h\nu^{-1})\mu
\end{equation}
for all $g,h\in Q$. Then the map $\tau\mapsto \tau^{(\nu,\mu)}$ is an action of $\operatorname{Comp}(Q,M)$ on $Z^2(Q,M)$
which preserves $B^2(Q,M)$ and so yields an action on $H^2(Q,M)$.

A {\em $KQ$-module extension} on $M$ by $Q$ is a group $E$  that fits in the short exact sequence

\begin{equation}\label{kge}
0\to M \stackrel{\iota}{\to} E \stackrel{\pi}{\to} Q \to 1
\end{equation}
so that the conjugation of $M$ (identified with $M\iota$) by elements of $E$ agrees with the $KQ$-module structure of $M$,
i.\,e. $m^e=m(e\pi)$ for all $m\in M$, $e\in E$.

\begin{prop}\cite[\S 2.7.4]{05HolEicObr} \label{comp} The classes of those isomorphisms of $KQ$-module extensions of $M$ by $Q$
that leave $M$ invariant as a $K$-module are in
a one-to-one correspondence with the orbits of $\operatorname{Comp}(Q,M)$ on $H^2(Q,M)$.
\end{prop}

In Proposition \ref{comp}, an isomorphism leaving $M$ invariant as a $K$-module means one that induces on $M$ an element of $\operatorname{Aut}_K(M)$.

The $KQ$-module structure on $M$ gives rise to the representation homomorphism $\mathcal{C}:Q\to \operatorname{Aut}_K(M)$ by the rule $\mathcal{C}(g):m\mapsto mg$ for all $m\in M$, $g\in Q$. Let $C$ be the centraliser of $\mathcal{C}(Q)$ in $\operatorname{Aut}_K(M)$. Then $(1,\gamma)\in \operatorname{Comp}(Q,M)$ for every $\gamma \in C$, because
$$(mg)\gamma=m\mathcal{C}(g)\gamma = m\gamma\mathcal{C}(g)=(m\gamma)g$$
for all $m\in M$, $g\in Q$. Hence, we also have an action of $C$ on both $Z^2(Q,M)$ and $H^2(Q,M)$ by setting
$\tau^\gamma=\tau^{(1,\gamma)}$ for $\tau \in Z^2(Q,M)$, $\gamma \in C$, i.\,e. $\tau^\gamma(g,h)=\tau(g,h)\gamma$. By Proposition \ref{comp}, this yields the following:

\begin{lem} \label{h2gscal} The elements of $H^2(Q,M)$ that are in the same $C$-orbit correspond to isomorphic
$KQ$-module extensions.
\end{lem}

In particular, we have the following fact, where two elements of $H^2(Q,M)$ are called scalar multiples if they differ by a factor in~$K^\times$.

\begin{cor} \label{h2scal} $KQ$-module extensions of $M$ by $Q$ corresponding to scalar multiples in $H^2(Q,M)$
are isomorphic.
\end{cor}

\section{Automorphisms of extensions}\label{sec:ae}

Fix an extension
\begin{equation}\label{ext}
 \bm{e}: \qquad 0\to M \stackrel{\iota}{\longrightarrow}E \to Q \to 1
\end{equation}
with abelian kernel $M$. Let $\mathcal{C}:Q\to \operatorname{Aut}(M)$ be the induced representation and
let $\overline{\varphi}\in H^2(Q,M)$ be the element that corresponds to $\bm{e}$.
We assume that $\mathcal{C}$ is faithful.
In particular, $Q\cong \mathcal{C}(Q)$ and the conjugation of $\mathcal{C}(Q)$ by any $\mu\in N_{\operatorname{Aut}(M)}(\mathcal{C}(Q))$
induces an element $\mu'\in \operatorname{Aut}(Q)$, i.\,e. $\mathcal{C}(g)^\mu=\mathcal{C}(g\mu')$ for
all $g\in Q$.
One defines an action of $N_{\operatorname{Aut}(M)}(\mathcal{C}(Q))$ on $H^2(Q,M)$ given by
\begin{equation}\label{act}
\overline{\psi}\mapsto(\mu')^{-1}\overline{\psi}\mu
\end{equation}
for every $\mu \in N_{\operatorname{Aut}(M)}(\mathcal{C}(Q))$ and $\overline{\psi}\in H^2(Q,M)$,
which should be understood modulo $B^2(Q,M)$ for representative cocycles, see \cite{82Rob} for details.
We denote by
$N_{\operatorname{Aut}(M)}^{\,\overline{\varphi}}(\mathcal{C}(Q))$ the stabiliser of
$\overline{\varphi}$ with respect to this action.
Let $\operatorname{Aut}(\bm{e})$ denote
the group of those automorphisms of $E$ that leave $M\iota$ invariant as a set.

\begin{prop}\cite[Statements (4.4),(4.5)]{82Rob}\label{Rob}
Let the extension $(\ref{ext})$ have an abelian kernel $M$ and let it
determine an element $\overline{\varphi}\in H^2(Q,M)$
and an injective induced representation $\mathcal{C}:Q\to \operatorname{Aut}(M)$.
Then there exists a short exact sequence of groups
\begin{equation}\label{ZAN}
0\to Z^1(Q,M)\to \operatorname{Aut}(\bm{e})\to N_{\operatorname{Aut}(M)}^{\,\overline{\varphi}}(\mathcal{C}(Q))\to 1.
\end{equation}
\end{prop}

\noindent
{\em Remark.} It is easy to see that, in the notation above, there is an embedding $N_{\operatorname{Aut}(M)}(\mathcal{C}(Q))\to \operatorname{Comp}(Q,M)$, $\mu\mapsto (\mu',\mu)$, where
we view $M$ as a $\mathbb{Z}Q$-module, under which action (\ref{act}) becomes a particular case of (\ref{cact}), and that this embedding is in fact an isomorphism in case $\mathcal{C}$ is faithful (which we assume).

\section{Cohomology of $PSL_2(q)$ in characteristic dividing~$q-1$}

The aim of this section is to classify up to group isomorphism nonsplit ex\-ten\-sions~(\ref{g}),
where $L=\operatorname{PSL}_2(q)$, $V$ is an elementary abelian $p$-group with irreducible induced action of $L$,
and $p\ne 2$ is a divisor of $q-1$.

By Lemma \ref{lgag}, $V$ must belong to the the principal $p$-block of $L$. This block
contains only one nonprincipal module with Brauer character $\chi$, see \cite{76Bur}.
The values of characters in the principal block are shown in Table \ref{bt}.

\begin{table}[htb]
\centering
\caption{Brauer $p$-modular characters of $L=PSL_2(q)$ in the principal block, where $2\ne p\mid(q-1)$.
Notation:  $q=l^m$, $l$ prime, $d=(2,q-1)$,
$x,y\in L$, $|x|=\frac{1}{d}(q-1)_{p'}$, $|y|=\frac{1}{d}(q+1)$.\label{bt}}
\begin{tabular}{c|rrrrrr}
            $q$ odd                &  $1a$      & $2a$                    & $la$ & $lb$ & $(x^r)^L$ & $(y^t)^L$  \\
     \hline
    $1_L^{\vphantom{A^A}}$    & $1$      & $ 1$    &  $1$  &  $1$ & $1$ & $1$ \\
    $\chi$                    & $q $     & $-1$    &  $0$    & $0$ & $1$  &  $-1$
\end{tabular}
\qquad
\begin{tabular}{c|rrrr}
            $q$ even                &  $1a$      & $2a$   & $(x^r)^L$ & $(y^t)^L$  \\
     \hline
    $1_L^{\vphantom{A^A}}$    & $1$      & $ 1$    & $1$ & $1$ \\
    $\chi$                    & $q $     & $0$     & $1$  &  $-1$
\end{tabular}
\end{table}

We first note that $V$ is not the principal module. Indeed, extension (\ref{g}) would otherwise be central,
but $\operatorname{Sch}(L)$ has no $p$-torsion, because
\begin{equation}\label{Sch}
\operatorname{Sch}(L)=\left\{
             \begin{array}{rl}
               \mathbb{Z}_2, & q\ne 9\  \text{odd or}\ q=4;\\
               \mathbb{Z}_6, & q=9;\\
               1, & q\ne 4 \ \text{even}
             \end{array}
           \right.
\end{equation}
as follows from \cite{85Atlas}. Therefore, $V$ must be the
$\mathbb{F}_pL$-module with character $\chi$.

We can now prove Theorem \ref{2coh} stated in the introduction.

\begin{proof} Let $P$ be the permutation $\mathbb{F}_pL$-module of dimension $q+1$ that
corresponds to the natural permutation action of $L$ on the projective line over $\mathbb{F}_q$.
We have $P=I_L\oplus V$, where $I_L$ is
the principal $\mathbb{F}_pL$-module. This can be deduced either by
considering the Brauer character $\chi$ of $V$ or from \cite[Table 1]{78Mor}. In particular,
by Lemma \ref{esum}, we have
\begin{equation}\label{hadd}
H^i(L,P)\cong H^i(L,I_L) \oplus H^i(L,V)
    \end{equation}
for $i=1,2$, since $H^i(L,B)\cong \operatorname{Ext}^i_{\mathbb{F}_pL}(\mathbb{F}_p,B)$
for every $\mathbb{F}_pL$-module $B$, see \cite[Exercise 6.1.2]{94Wei}.
Since $P$ is a permutation module, we have $P\cong (I_H)^L$, where $I_H$ is the principal
$\mathbb{F}_pH$-module for a point stabiliser $H\leqslant L$. Hence, Lemma \ref{shap} implies
\begin{equation}\label{shs}
H^i(L,P)\cong H^i(H,I_H)
\end{equation}
for $i=1,2$. By Lemma \ref{coh1}$(i)$, we have
\begin{equation}\label{h1l}
H^1(L,I_L)\cong \operatorname{Hom}(L/L',I_L)=0,
\end{equation}
since $L=L'$. Also,
\begin{equation}\label{h1h}
H^1(H,I_H)\cong \operatorname{Hom}(H/H',I_H)\cong \mathbb{F}_p,
\end{equation}
since $I_H\cong \mathbb{F}_p$, $H\cong \mathbb{F}_q \leftthreetimes \mathbb{Z}_{(q-1)/(2,q-1)}$, and $p\mid (q-1)$. By Lemma \ref{ucoe}, we have
$$
H^2(L,I_L)\cong \operatorname{Hom}(H_2(L,\mathbb{Z}),I_L)\oplus \operatorname{Ext}(H_1(L,\mathbb{Z}),I_L),
$$
where the first summand vanishes, since $H_2(L,\mathbb{Z})\cong \operatorname{Sch}(L)$ has no $p$-torsion by~(\ref{Sch}), and the second summand vanishes by Lemma \ref{coh1}$(ii)$, since $L/L'=1$. Thus
\begin{equation}\label{h2lil}
H^2(L,I_L)=0.
\end{equation}
Finally, Lemma \ref{ucoe} also yields
\begin{equation}\label{h2hihe}
H^2(H,I_H)\cong \operatorname{Hom}(H_2(H,\mathbb{Z}),I_H)\oplus \operatorname{Ext}(H_1(H,\mathbb{Z}),I_H).
\end{equation}
By Lemma \ref{tim}, the $p$-part of $H_2(H,\mathbb{Z})\cong \operatorname{Sch}(H)$ is isomorphic to a subgroup of $\operatorname{Sch}(S)$ for a $p$-Sylow subgroup $S$ of $H$. However, $S$ is cyclic and
cyclic groups have trivial Schur multiplier. Thus, the first summand in (\ref{h2hihe}) vanishes, because $I_H\cong \mathbb{F}_p$. Since $H_1(H,\mathbb{Z})\cong H/H'\cong \mathbb{Z}_{(q-1)/d}$ and
$\operatorname{Ext}(\mathbb{Z}_{(q-1)/d},\mathbb{F}_p)\cong \mathbb{F}_p$ by Lemma \ref{ezmn}, we have
\begin{equation}\label{h2hih}
H^2(H,I_H)\cong \mathbb{F}_p.
\end{equation}
The claim follows by combining (\ref{hadd}) through (\ref{h2hih}).
\end{proof}

We can now prove Theorem \ref{nonsplit} stated in the introduction.

\begin{proof} As we explained in the beginning of this section, $V$ viewed as an $\mathbb{F}_pL$-module
must be the unique nonprincipal module in the principal $p$-block of $L$. This module has dimension $q$
and can be written over $\mathbb{F}_p$, since it is a direct summand of a permutation module. Therefore,
$|V|=p^q$. By Theorem \ref{2coh}, we have $H^2(V,L)\cong \mathbb{F}_p$ and so all nonzero elements of
$H^2(V,L)$ are scalar multiples of one another. By Corollary \ref{h2scal}, they correspond to isomorphic nonsplit extensions. The claim follows.
\end{proof}

\section{The automorphism group}

In this section, we prove that the structure of the automorphism group of the unique nonsplit
extension from Theorem \ref{nonsplit} is as stated in Theorem \ref{aut}.

\begin{proof}
Consider the extension $\bm{e}$ given by (\ref{g}).  Theorem \ref{nonsplit} implies that
$G$ is unique up to isomorphism and $V$ has order $p^q$. Moreover, viewed as an $\mathbb{F}_p L$-module, $V$ has Brauer character $\chi$ from Table \ref{bt}.
By Proposition \ref{Rob}, we have the short exact sequence
\begin{equation}\label{ZAN1}
0\to Z^1(L,V)\to \operatorname{Aut}(\bm{e})\to N_{\operatorname{Aut}(V)}^{\,\overline{\varphi}}(\mathcal{X}(L))\to 1,
\end{equation}
where the representation $\mathcal{X}:L\to \operatorname{Aut}(V)$
and the element $\overline{\varphi}\in H^2(L,V)$ are determined by
(\ref{g}). First, note that $\operatorname{Aut}(\bm{e})=\operatorname{Aut}(G)$ as $V$ is characteristic in~$G$. Denote $W=Z^1(L,V)$. Since $B^1(L,V)\cong V/C_V(L)$ and $L$ acts on $V$ irreducibly and nontrivially, we have $C_V(L)=0$ and $B^1(L,V)\cong V$. Now, since $H^1(L,V)=Z^1(L,V)/B^1(L,V)$, we have
$|Z^1(L,V)|=p^{q+1}$ in view of Theorem \ref{2coh}.

Denote $N=N_{\operatorname{GL}(V)}(\mathcal{X}(L))$ and $Z=C_{\operatorname{GL}(V)}(\mathcal{X}(L))$.
By Proposition \ref{normim}, we have $N/Z\cong I_{\operatorname{Aut}(L)}(\chi)$.
Since $\chi$ is the only irreducible character of $L$ of dimension~$q$,
it must be invariant under any automorphism; in particular,
$I_{\operatorname{Aut}(L)}(\chi)=\operatorname{Aut}(L)$.
By \cite{85Atlas},  $\operatorname{Aut}(L)\cong P\Gamma L_2(q)$.
Since $V$ is absolutely irreducible as an $\mathbb{F}_pL$-module, by Schur's lemma, we see that
$Z\cong \mathbb{F}_p^\times\cong \mathbb{Z}_{p-1}$ consists of scalars.

In order to determine the structure of the stabiliser $N_0=N_{\operatorname{Aut}(V)}^{\,\overline{\varphi}}(\mathcal{X}(L))$,
we consider the action of $N$ on $H^2(L,V)$ as explained in Section \ref{sec:ae}. Let $H^\times$ denote the set
of $p-1$ nonzero elements of $H^2(L,V)$. The elements of $H^\times$ correspond to nonsplit extensions and so we have an action
homomorphism $\alpha: N\to \operatorname{Sym}(H^\times)$ to the symmetric group on $H^\times$.
Since all nonsplit extension of $V$ by $L$
are isomorphic by Theorem~\ref{nonsplit}, we may assume that $\overline{\varphi}$ is an arbitrary element of $H^\times$.
The subgroup $Z\leqslant N$ acts on $H^\times$ by
scalar multiplication, cf. Corollary \ref{h2scal}, and so the image $\alpha(Z)$ is a cyclic subgroup
of $\operatorname{Sym}(H^\times)$ generated by a full cycle of length $p-1$.
Since $Z$ is central in $N$, $\alpha(N)$ must centralise $\alpha(Z)$. However, a full cyclic subgroup is self-centralising in
$\operatorname{Sym}(H^\times)$ and so $\alpha(Z)$ must be the entire image $\alpha(N)$. Thus, $\operatorname{Ker}(\alpha)$ is a normal
subgroup of $N$ of index $p-1$ which intersects trivially with $Z$ and is thus isomorphic to $N/Z\cong P\Gamma L_2(q)$.
Furthermore, $\operatorname{Ker}(\alpha)$ coincides with the stabiliser of every element of $H^\times$ which yields
$N=N_0\times Z$ and $N_0\cong P\Gamma L_2(q)$ as claimed.
\end{proof}

It also follows from this proof that the representation $\mathcal{X}:L\to \operatorname{Aut}(V)$
with character $\chi$ extends to a representation of $I_{\operatorname{Aut}(L)}(\chi)\cong P\Gamma L_2(q)$. This fact does
not hold in general for a simple group $L$ and its irreducible character $\chi$, see \cite[Example 1]{21RevZav}.

{\em Acknowledgement.\/} The work was supported by the RAS Fundamental Research Program, project FWNF-2022-0002.

\end{document}